\documentclass[useAMS]{my_gCOM2e}

\usepackage{url}

\begin{document}

\markboth{H.S. Rodrigues, M.T.T. Monteiro, D.F.M. Torres}{International Journal of Computer Mathematics}

\title{Bioeconomic Perspectives to an Optimal Control Dengue Model}

\author{Helena Sofia Rodrigues$^{\rm a}$, M. Teresa T. Monteiro$^{\rm b}$$^{\ast}$\thanks{$^\ast$Corresponding author.
Email: tm@dps.uminho.pt\vspace{6pt}}, Delfim F. M. Torres$^{\rm c}$\\\vspace{6pt}
$^{\rm a}${\em{School of Business Studies,
Viana do Castelo Polytechnic Institute, Portugal}}; $^{\rm b}${\em{Algoritmi Centre,
Department of Production and Systems, University of Minho, Braga, Portugal}};
$^{\rm c}${\em{Center for Research and Development in Mathematics and Applications,
Department of Mathematics, University of Aveiro, Aveiro, Portugal}}\\\vspace{6pt}
\received{Submitted 09-Aug-2012; revised 07-Jan-2013; accepted 21-Mar-2013}}

\maketitle


\begin{abstract}
A model with six mutually-exclusive compartments related to dengue is studied.
Three vector control tools are considered: insecticides (larvicide and adulticide)
and mechanical control. The basic reproduction number associated to the model is presented.
The problem is studied using an optimal control approach.
The human data is based on the dengue outbreak that occurred in Cape Verde.
Control measures are simulated in different scenarios and their consequences analyzed.

\bigskip

\begin{keywords}
modeling; optimal control; basic reproduction number; vector control; dengue
\end{keywords}

\begin{classcode}
49M37; 92B05
\end{classcode}

\end{abstract}


\section{Introduction}

Dengue is a vector borne disease transmitted to humans by the bite
of an infected female \emph{Aedes} mosquito. Dengue transcends
international borders and can be found in tropical and sub-tropical
regions around the world, predominantly in urban and semi-urban areas.
The risk may be aggravated further due to climate changes and
globalization, as a consequence of the huge volume of international
tourism and trade \cite{Semenza2009}.

There are four distinct, but closely related, viruses that cause
dengue. Recovery from infection by one virus provides lifelong
immunity against that virus but confers only partial and transient
protection against subsequent infection by the other three
viruses \cite{Wearing2006}. Unfortunately,
there is no specific effective treatment for dengue.

Primary prevention of dengue resides mainly in mosquito control.
There are two main methods: larval control and adult mosquito
control, depending on the intended target \cite{Natal2002}.
The application of adulticides can have a powerful impact
on the abundance of adult mosquito vector.
This is the most common measure. However, the efficacy
is often constrained by the difficulty in achieving sufficiently
high coverage of resting surfaces. Besides,
the long term use of adulticide has several risks:
the resistance of the mosquito to the product, reducing its efficacy,
and the killing of other species that live in the same habitat.
Larvicide treatment is done through a long-lasting chemical, in order
to kill larvae, preferably with WHO clearance for use
in drinking water \cite{Derouich2003}. Larvicide treatment
is an effective way to control the vector larvae, together
with \emph{mechanical control}, which is related with
educational campaigns. The mechanical control must be done both
by public health officials and by residents in affected areas.
The participation of the entire population is essential
to remove still water from domestic recipients,
eliminating possible breeding sites \cite{Who2009}.
The SIR+ASI model considered here was studied 
in a previous paper \cite{Sofia2013}, but only using the ODE system:
an optimal control approach is here considered for the first time.
This provides a new different mathematical perspective to the subject.
Furthermore, a numerical procedure, varying the control weights in the model, is performed here
in order to evaluate the control that is most effective in the design of optimal strategies.

The paper is organized as follows. Next Section~\ref{sec:2}
presents the mathematical model under study.
Section~\ref{sec:3} is concerned with the basic reproduction number of the model.
The optimal control approach is addressed in Section~\ref{sec:4} while
the computational experiments, using different situations, are reported
in Section~\ref{sec:5}. Finally, some conclusions are carried out
in Section~\ref{sec:6}.


\section{Mathematical Model}
\label{sec:2}

Mathematical modeling is an interesting tool for
understanding epidemiological diseases and for proposing
effective strategies to fight them \cite{Lenhart2007}.
Taking into account the model presented in \cite{Dumont2010,Dumont2008}
for the chikungunya disease, and the considerations of \cite{Sofia2009, Sofia2010c},
a new model, more adapted to the dengue reality, is proposed.
Our epidemiological model for dengue is similar to models for chikungunya
because the vector is from the same mosquito family: the family of \emph{Aedes} mosquitoes.
In chikungunya the main vector is \emph{Aedes albopictus} while in dengue
the vector is \emph{Aedes aegypti}. However, here some parameter fitting is done, taking into
account the human and the vector populations and their specific
features. More precisely, we use data from Cape Verde, while in \cite{Dumont2010,Dumont2008}
the data is from Reunion islands. Moreover, the models used in \cite{Dumont2010,Dumont2008}
only consider the ODE system. In contrast, here we use the optimal control approach
in the epidemiological model.

The notation used in our mathematical model includes three
epidemiological states for humans, indexed by $h$:
\begin{equation*}
\begin{split}
S_h & \text{ --- susceptible (individuals who can contract the disease); }\\
I_h & \text{ --- infected (individuals capable of transmitting the disease to others); }\\
R_h & \text{ --- resistent (individuals who have acquired immunity). }
\end{split}
\end{equation*}

It is assumed that the total human population $N_h$
is constant along time: $N_h=S_h(t)+I_h(t)+R_h(t)$.

There are three other state variables, related
to the female mosquitoes, indexed by $m$:
\begin{equation*}
\begin{split}
A_m & \text{ --- aquatic phase (that includes the egg, larva and pupa stages);}\\
S_m & \text{ --- susceptible (mosquitoes that are able to contract the disease);}\\
I_m & \text{ --- infected (mosquitoes capable of transmitting the disease to humans).}
\end{split}
\end{equation*}

Due to short lifetime of mosquitoes (approximately 10 days), there is no resistant phase.
Humans and mosquitoes are assumed to be born susceptible.

To analyze the effect of campaigns in the combat
of the mosquito, three controls are considered:\footnote{The control $\alpha$ cannot be zero 
because it appears in the denominator of a fraction in the ODE system \eqref{cap6_ode1}.}
\begin{equation*}
\begin{split}
c_{A} & \text{ --- proportion of larvicide ($0\leq c_A\leq1$);}\\
c_{m} & \text{ ---  proportion of adulticide ($0\leq c_m\leq1$);}\\
\alpha & \text{ ---  proportion of mechanical control
($0 < \alpha \leq1$).}
\end{split}
\end{equation*}

The aim of this work is to simulate different realities in order
to find the best policy to decrease the number of infected human.
A temporal mathematical model is introduced, with mutually-exclusive compartments,
to study the outbreak of 2009 in Cape Verde islands and
improving the model described in \cite{Sofia2009}.
The model considers the following parameters:
\begin{equation*}
\begin{split}
N_h & \text{ --- total human population;}\\
B & \text{ --- average daily biting (per day);}\\
\beta_{mh} & \text{ --- transmission probability from $I_m$ (per bite);}\\
\beta_{hm} & \text{ --- transmission probability from $I_h$ (per bite);}\\
1/\mu_{h} & \text{ --- average lifespan of humans (in days);}\\
1/\eta_{h} & \text{ --- mean viremic period (in days);}\\
1/\mu_{m} & \text{ --- average lifespan of adult mosquitoes (in days);}\\
\varphi & \text{ --- number of eggs at each deposit per capita (per day);}\\
1/\mu_{A} & \text{ --- natural mortality of larvae (per day);}\\
\eta_{A} & \text{ --- maturation rate from larvae to adult (per day);}\\
m & \text{ --- female mosquitoes per human;}\\
k & \text{ --- number of larvae per human.}
\end{split}
\end{equation*}

The dengue epidemic is modeled by the following time-varying
nonlinear system of differential equations:
\begin{equation}
\label{cap6_ode1}
\left\{
\begin{split}
\frac{dS_h}{dt}(t) &= \mu_h N_h - \left(B\beta_{mh}\frac{I_m}{N_h}+\mu_h\right)S_h\\
\frac{dI_h}{dt}(t) &= B\beta_{mh}\frac{I_m}{N_h}S_h -(\eta_h+\mu_h) I_h\\
\frac{dR_h}{dt}(t) &= \eta_h I_h - \mu_h R_h\\
\frac{dA_m}{dt}(t) &= \varphi \left(1-\frac{A_m}{\alpha k N_h}\right)(S_m+I_m)
-\left(\eta_A+\mu_A + c_A\right) A_m\\
\frac{dS_m}{dt}(t) &= \eta_A A_m
-\left(B \beta_{hm}\frac{I_h}{N_h}+\mu_m + c_m\right) S_m\\
\frac{dI_m}{dt}(t) &= B \beta_{hm}\frac{I_h}{N_h}S_m
-\left(\mu_m + c_m\right) I_m
\end{split}
\right.
\end{equation}
with the initial conditions
\begin{equation}
\label{cap6_initial}
\begin{tabular}{llll}
$S_h(0)=S_{h0},$ &  $I_h(0)=I_{h0},$ &
$R_h(0)=R_{h0},$ \\
$A_m(0)=A_{m0},$ & $S_{m}(0)=S_{m0},$ & $I_m(0)=I_{m0}.$
\end{tabular}
\end{equation}


\section{Basic Reproduction Number}
\label{sec:3}

An important measure of transmissibility of the disease
is given by the basic reproduction number. It represents the
expected number of secondary cases produced in a completed
susceptible population, by a typical infected individual during
its entire period of infectiousness \cite{Hethcote2000}.
It can be shown (see \cite{Sofia2013}) that
the basic reproduction number $\mathcal{R}_0$ associated
to \eqref{cap6_ode1} is given by
\begin{equation}
\label{eq:R0}
\mathcal{R}_0 = \left(\frac{\alpha k B^2 \beta_{hm} \beta_{mh}
\mathcal{M}}{\varphi (\eta_h + \mu_h) (c_m + \mu_m)^2}\right)^{\frac{1}{2}}.
\end{equation}
The model has two different populations (host and vector)
and the expected basic reproduction number reflects the
infection transmitted from host to vector and vice-versa.
If $\mathcal{R}_0<1$, then, on average, an infected individual
produces less than one new infected individual over the course of
its infectious period, and the disease cannot grow. Conversely, if
$\mathcal{R}_0>1$, then each individual infects more than one
person, and the disease invades the population.

\begin{figure}
\centering
\subfigure[$\mathcal{R}_{0}$ as a function of $c_m$ and $c_A$]{\label{cap6_R0_cm_CA}
\includegraphics[width=0.5\textwidth]{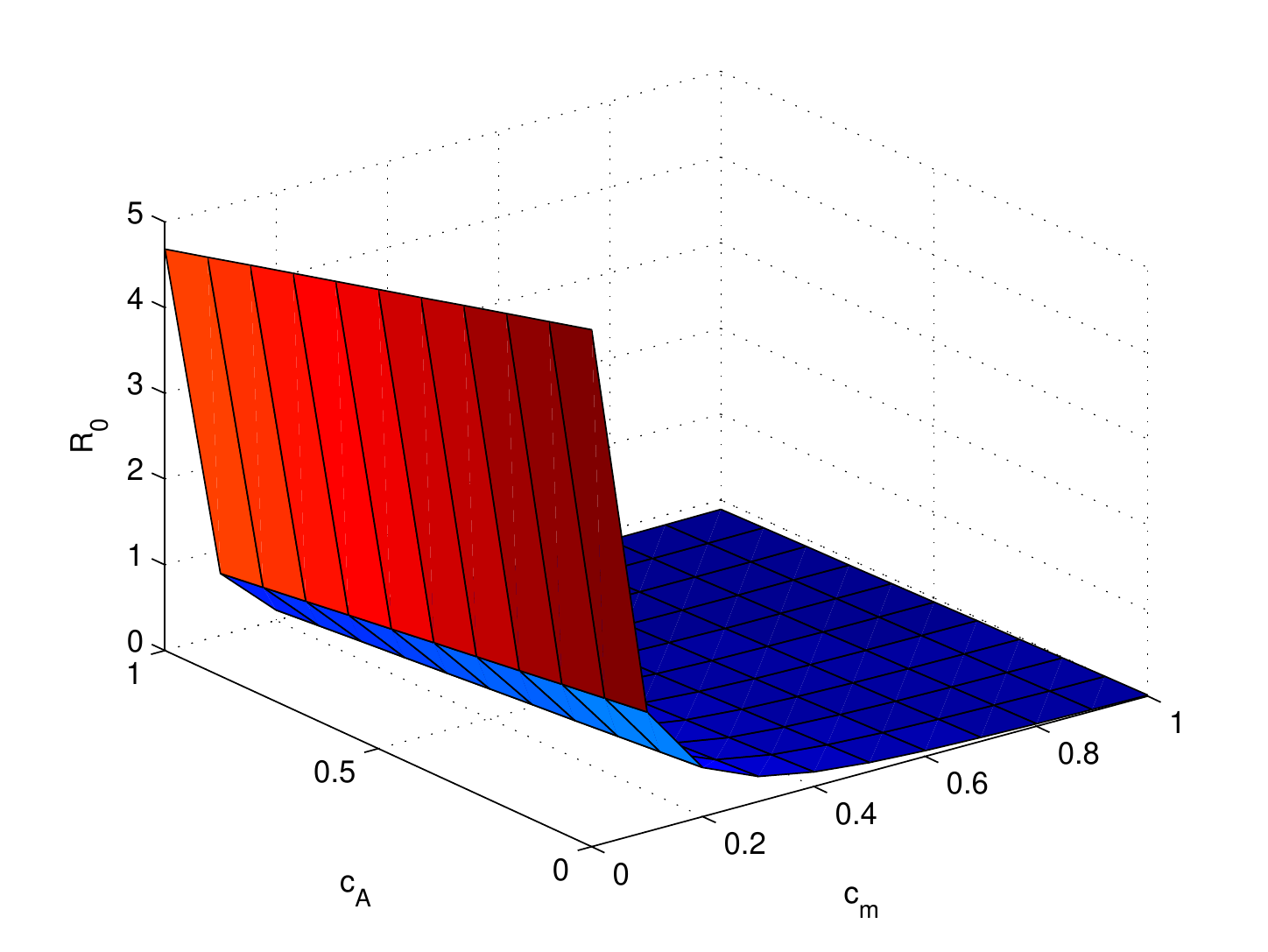}}\\
\subfigure[$\mathcal{R}_{0}$ as a function of $c_m$ and $\alpha$]{\label{cap6_R0_cm_alpha}
\includegraphics[width=0.5\textwidth]{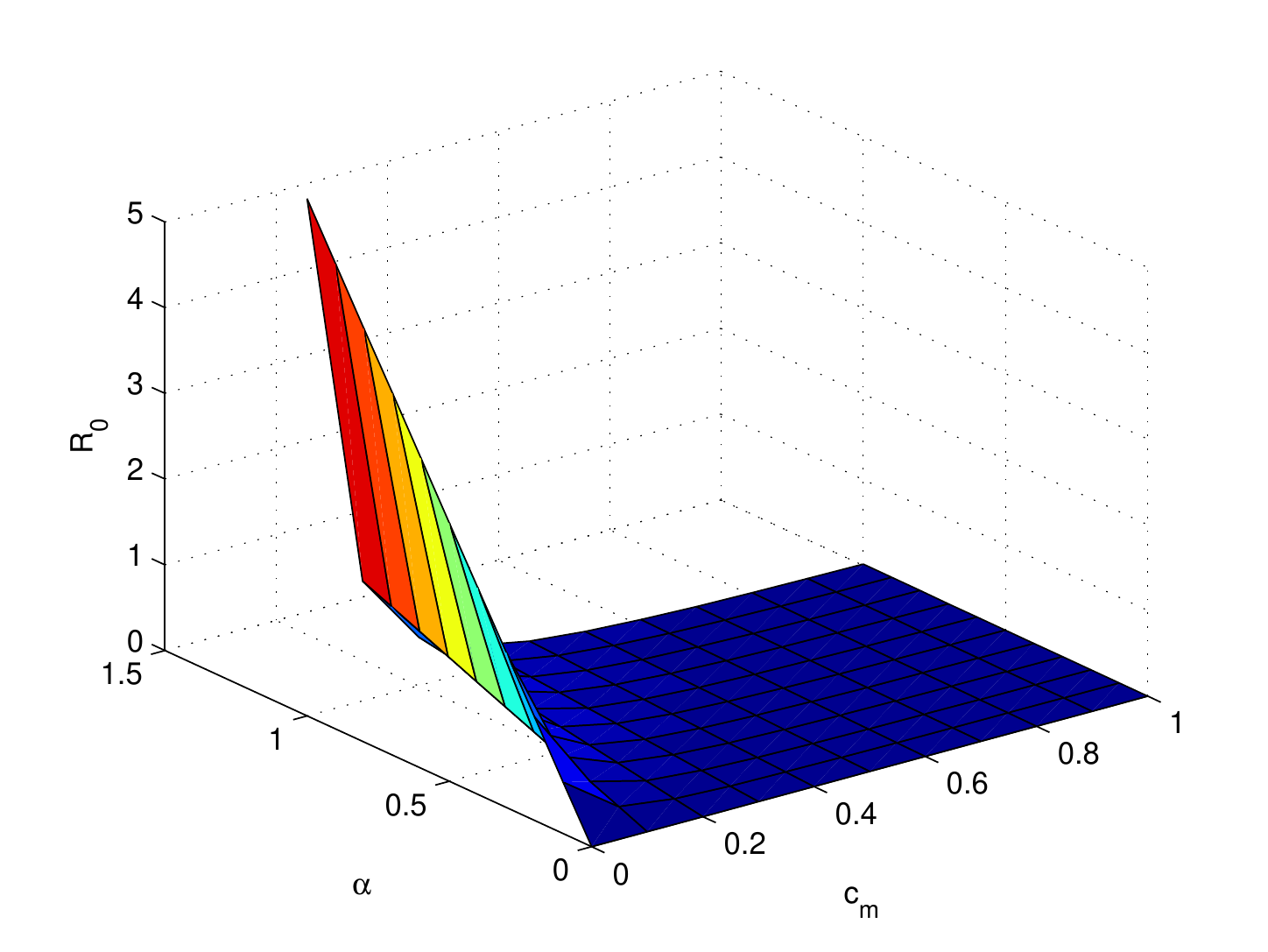}}\\
\subfigure[$\mathcal{R}_{0}$ as a function of $c_A$ and $\alpha$]{\label{cap6_R0_cA_alpha}
\includegraphics[width=0.5\textwidth]{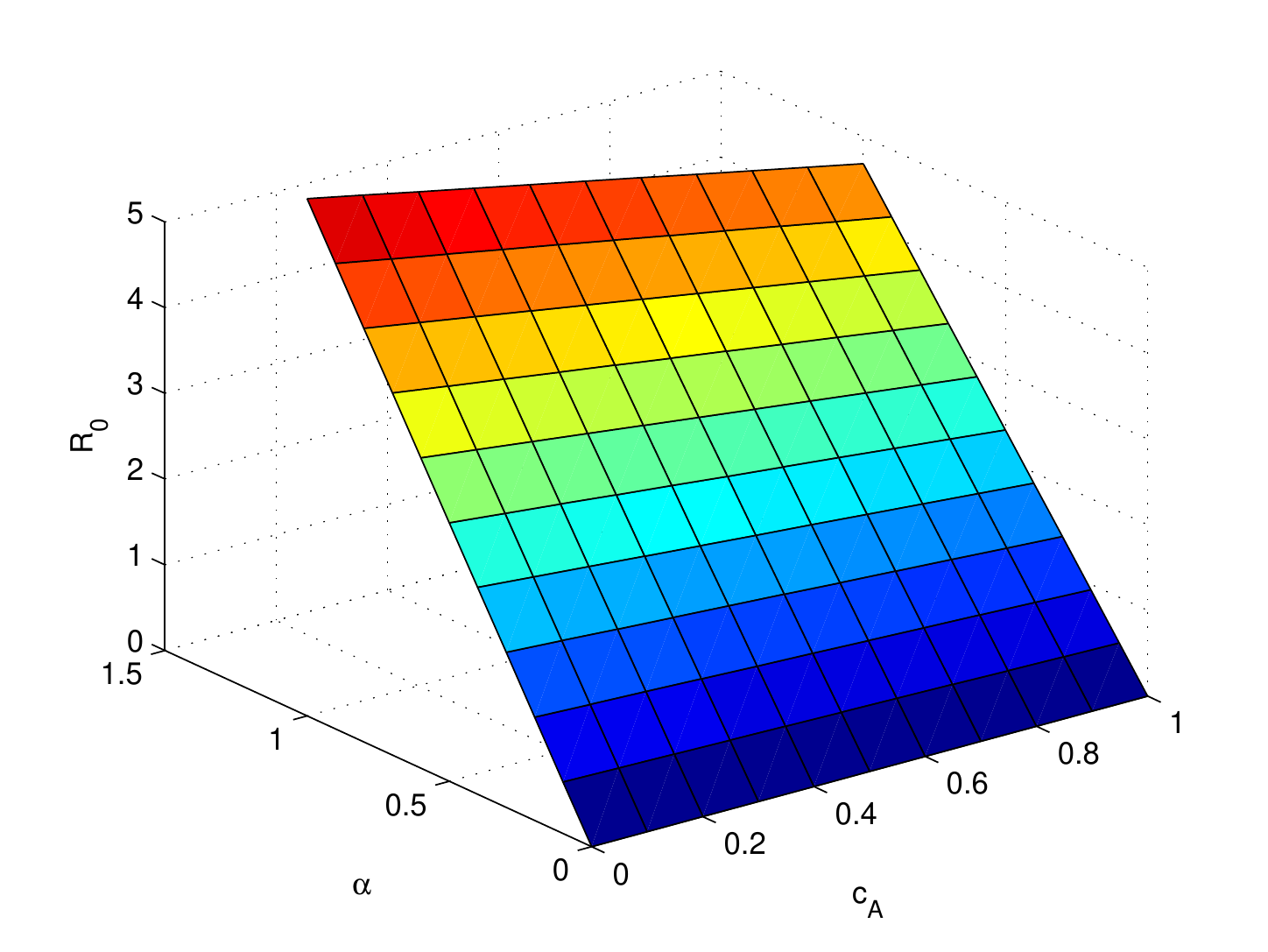}}
\caption{Influence of the controls on the basic reproduction number $\mathcal{R}_{0}$}
\label{cap6_R0_cm_CA_alpha}
\end{figure}

Assuming that parameters are fixed, the threshold $\mathcal{R}_0$
is influenceable by the control values. Figure~\ref{cap6_R0_cm_CA_alpha}
gives this relationship. We see that the control $c_m$
is the one that most influences the basic reproduction number to stay below unit.
Besides, the control in the aquatic phase alone is not enough
to maintain $\mathcal{R}_0$ below unit: it requires an application close to 100\%.


\section{Optimal Control Approach}
\label{sec:4}

Epidemiological models may give some basic guidelines for public health practitioners,
comparing the effectiveness of different potential management strategies.
In reality, a range of constraints and trade-offs may substantially influence
the choice of a practical strategy, and therefore their inclusion
in any modeling analysis may be important. Frequently, epidemiological models need
to be coupled to economic considerations, such that control strategies can be judged
through holistic cost-benefit analysis. Control of livestock disease is a scenario
when cost-benefit analysis can play a vital role in choosing between cheap,
weak controls that lead to a prolonged epidemic, or expensive
but more effective controls that lead to a shorter outbreak.
In our numerical simulations, the data from human initial conditions 
was obtained through the Ministry of Health from Cape Verde \cite{CapeVerdeSaude}. 
Since it was the first time that an outbreak of dengue occurred in Cape
Verde, there was no time to follow the mosquito evolution. Besides, the Health
authorities of Cape Verde believe that the mosquito came from Brazil, based on the intensive
commercial trade and migration between the two countries. Therefore,
the vector data was based on the Brazil reality \cite{Yang2009,Thome2010}.
Normalizing the previous ODE system \eqref{cap6_ode1}--\eqref{cap6_initial}, we obtain:
\begin{equation}
\label{cap6_ode_norm}
\left\{
\begin{split}
\dfrac{ds_h}{dt} &= \mu_h - \left(B\beta_{mh}m i_m+\mu_h\right)s_h\\
\dfrac{di_h}{dt} &= B\beta_{mh}m i_m s_h -(\eta_h+\mu_h) i_h\\
\dfrac{dr_h}{dt} &= \eta_h i_h - \mu_h r_h\\
\dfrac{da_m}{dt} &= \varphi \frac{m}{k}\left(1-\frac{a_m}{\alpha }\right)(s_m+i_m)
-\left(\eta_A+\mu_A + c_A\right) a_m\\
\dfrac{ds_m}{dt} &= \eta_A \frac{k}{m}a_m -\left(B \beta_{hm}i_h+\mu_m + c_m\right) s_m\\
\dfrac{di_m}{dt} &= B \beta_{hm}i_h s_m -\left(\mu_m + c_m\right) i_m
\end{split}
\right.
\end{equation}
with the initial conditions
\begin{equation}
\label{chap6_initial_norm}
\begin{tabular}{llll}
$s_h(0)=0.9999$, &  $i_h(0)=0.0001$, &
$r_h(0)=0$, \\
$a_m(0)=1$, & $s_{m}(0)=1$, & $i_m(0)=0$.
\end{tabular}
\end{equation}
A cost functional was introduced,
\begin{equation}
\label{chap6_functional}
J[c_A(\cdot),c_m(\cdot),\alpha(\cdot)]=\int_{0}^{t_f}\left[\gamma_D I_h(t)^2
+\gamma_S c_m(t)^2+\gamma_L c_A(t)^2+\gamma_E \left(1-\alpha\right)^2\right]dt,
\end{equation}
where $\gamma_D$, $\gamma_S$, $\gamma_L$ and $\gamma_E$ are weights related
to the costs of the disease, adulticide, larvicide and mechanical control, respectively.
In this way, an optimal control problem is defined:
\begin{equation*}
\begin{tabular}{ll}
minimize & \eqref{chap6_functional} \\
subject to & \eqref{cap6_ode_norm}, \quad \eqref{chap6_initial_norm}, \quad
$0\leq c_A \leq 1, \quad 0\leq c_m \leq 1, \quad 0< \alpha \leq 1$.
\end{tabular}
\end{equation*}


\section{Numerical Experiments with Three Controls}
\label{sec:5}

The simulations were carried out with the numerical values
$N_h=480000$, $B=0.8$, $\beta_{mh}=0.375$,
$\beta_{hm}=0.375$, $\mu_{h}=1/(71\times365)$, $\eta_{h}=1/3$,
$\mu_{m}=1/10$, $\varphi=6$, $\mu_A=1/4$, $\eta_A=0.08$, $m=3$,
$k=3$, and initial conditions \eqref{chap6_initial_norm}.
The optimal control problem was solved using two different approaches
and software packages: \textsf{DOTcvp} \cite{bib:DOTcvp}
and \textsf{Muscod-II} \cite{bib:Muscod-II}.
In both cases the simulations were similar. Thus,
only the \textsf{DOTcvp} results are reported here.
We remark that our results cannot be compared with those 
of \cite{Dumont2010,Dumont2008} for several reasons.
First, the parameters of both works are different because they describe two distinct
diseases and two different realities. Second, in \cite{Dumont2010,Dumont2008}
the optimal control approach is not considered and, therefore,
the optimal control strategy is not computed
(only simulations of some suboptimal strategies are done). The
study of the variation of control weights carried out here is,
due to economical constraints, a very important issue,
allowing the decision-makers (the health authorities) 
to define bioeconomic strategies.


\subsection{All controls with the same weight}

\begin{figure}[ptbh]
\begin{center}
\includegraphics[scale=0.7]{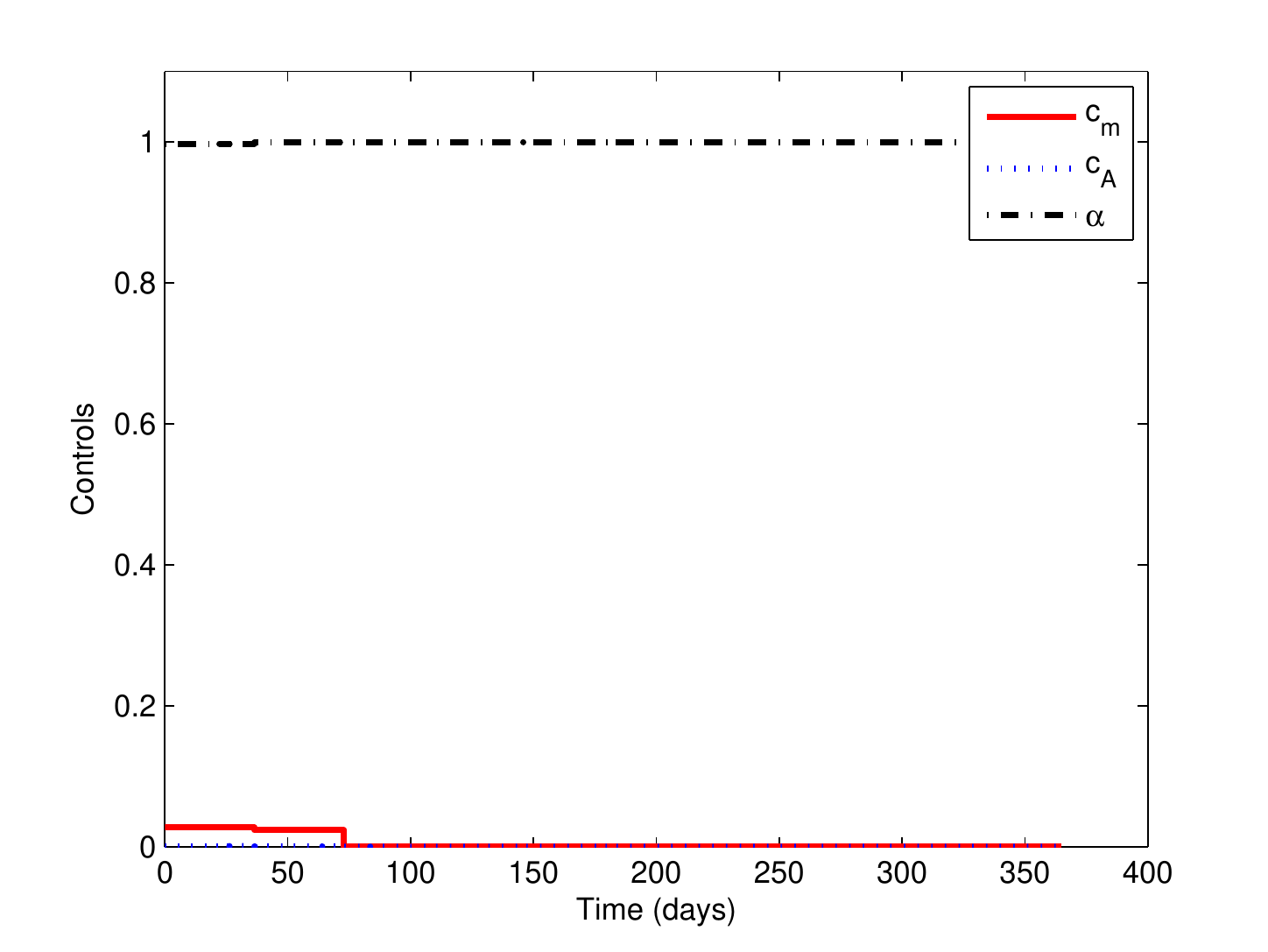}
\end{center}
\caption{Optimal control functions ($\gamma_D=\gamma_S=\gamma_L=\gamma_E=0.25$)}
\label{cap6_all_controls}
\end{figure}

To begin, the same weights were considered:
$\gamma_D= \gamma_S=\gamma_L=\gamma_E=0.25$.
The optimal functions for the controls are given
in Figure~\ref{cap6_all_controls}.

The adulticide was the control that more influences the decreasing
of the basic reproduction number \eqref{eq:R0} and, as a consequence,
the decreasing of the number of infected persons and mosquitoes.
Therefore, the adulticide was almost the one to be used.
The other controls do not assume here an important role in the epidemic episode,
because all the events happen in a short period of time,
which means that adulticide has more impact. However, the control
of the mosquito in the aquatic phase cannot be neglected.
In situations of longer epidemic episodes or even in an endemic situation,
the larval control represents an important tool.

\begin{figure}[ptbh]
\begin{center}
\includegraphics[scale=0.7]{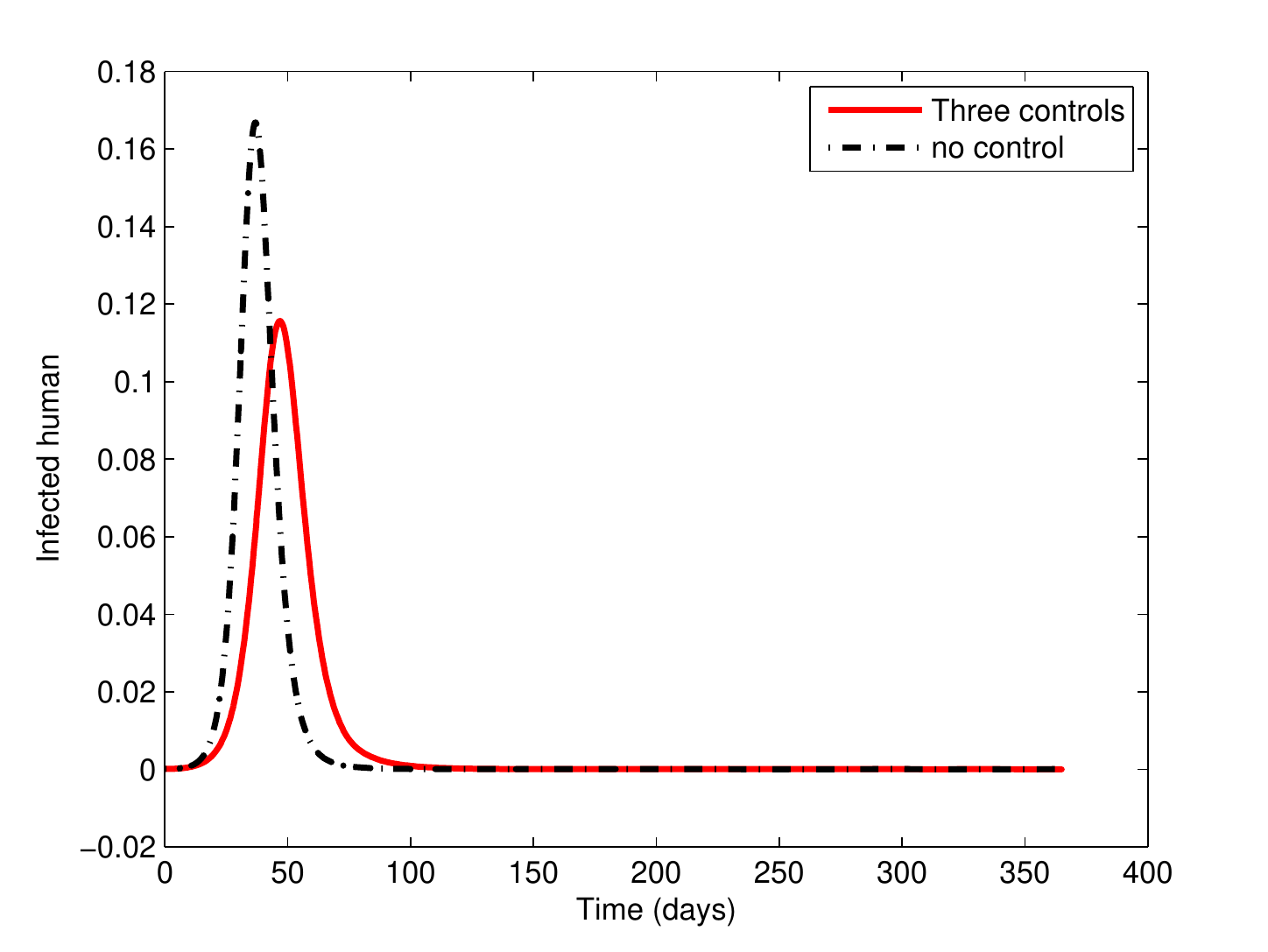}
\end{center}
\caption{Comparison of infected individuals under
an optimal control strategy with that of no controls.}
\label{cap6_all_controls_vs_nocontrol}
\end{figure}

Figure~\ref{cap6_all_controls_vs_nocontrol} presents the number of infected human.
Comparing the optimal control case with the situation of no controls,
the number of infected people decreased considerably. Besides, in the situation
where optimal control is used, the peak of infected people is minor,
which facilitates the work in health centers,
because they can provide a better medical monitoring.

\subsection{Controls with different weights}

A second analysis was made, taking into account
different weights on the functional \eqref{chap6_functional}.
Table~\ref{table_chap6_optimal weights} summarizes the weights chosen and
the associated perspectives. Not only economic issues
(cost of insecticides and educational campaigns) are considered,
but also human issues. In case A all costs are equal.
In case B more impact is given to the infected people, considering that the treatment
and absenteeism to work is very prejudicial to the country, when compared with the
cost of insecticides and educational campaigns. In case C, the costs of killing
mosquitoes and educational campaigns are the ones with more impact in the economy.

\begin{table}[ptbh]
\begin{center}
\small
\begin{tabular}{l l c}
\hline
& Values for weights & Cost obtained\\
\hline
Case A & $\gamma_D=0.25$; $\gamma_S=0.25$; $\gamma_L=0.25$; $\gamma_E=0.25$ & 0.06691425\\
Case B & $\gamma_D=0.55$; $\gamma_S=0.15$; $\gamma_L=0.15$; $\gamma_E=0.15$ & 0.10431186\\
Case C & $\gamma_D=0.10$; $\gamma_S=0.30$; $\gamma_L=0.30$; $\gamma_E=0.30$ & 0.03012849\\
\hline
\end{tabular}
\caption{Different weights for the functional \eqref{chap6_functional} and respective values}
\label{table_chap6_optimal weights}
\end{center}
\end{table}

The higher total costs were obtained when the human life has more weight
than the controls measures (Table~\ref{table_chap6_optimal weights}).

\begin{figure}[ptbh]
\begin{center}
\includegraphics[scale=0.7]{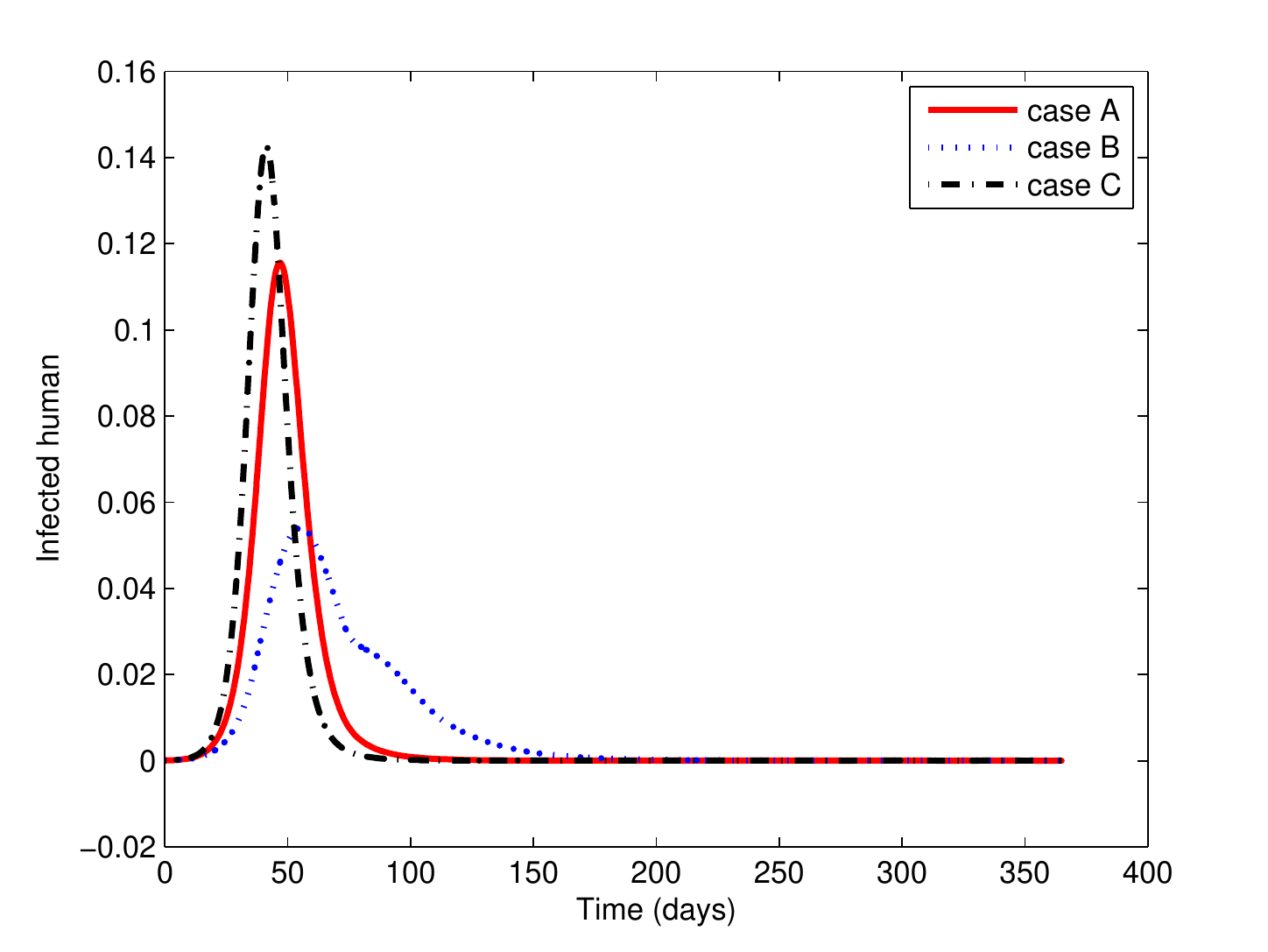}
\end{center}
\caption{Infected individuals in the three bioeconomic perspectives}
\label{cap6_human_infected_CasesABC}
\end{figure}

\begin{figure}
\centering
\subfigure[Adulticide]{\label{cap6_proportion_adulticide_CasesABC}
\includegraphics[width=0.55\textwidth]{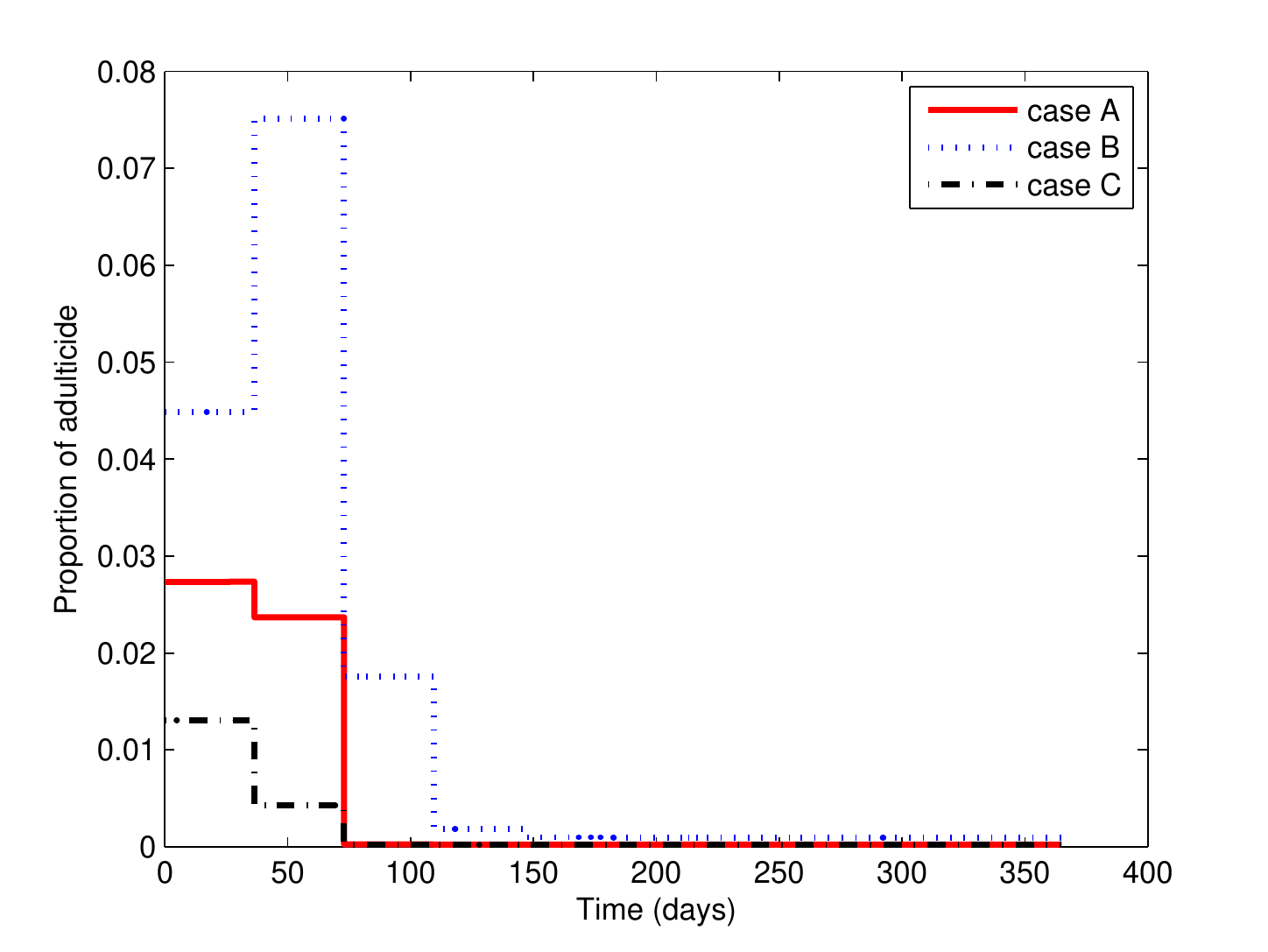}}\\
\subfigure[Larvicide]{\label{cap6_proportion_larvicide_CasesABC}
\includegraphics[width=0.55\textwidth]{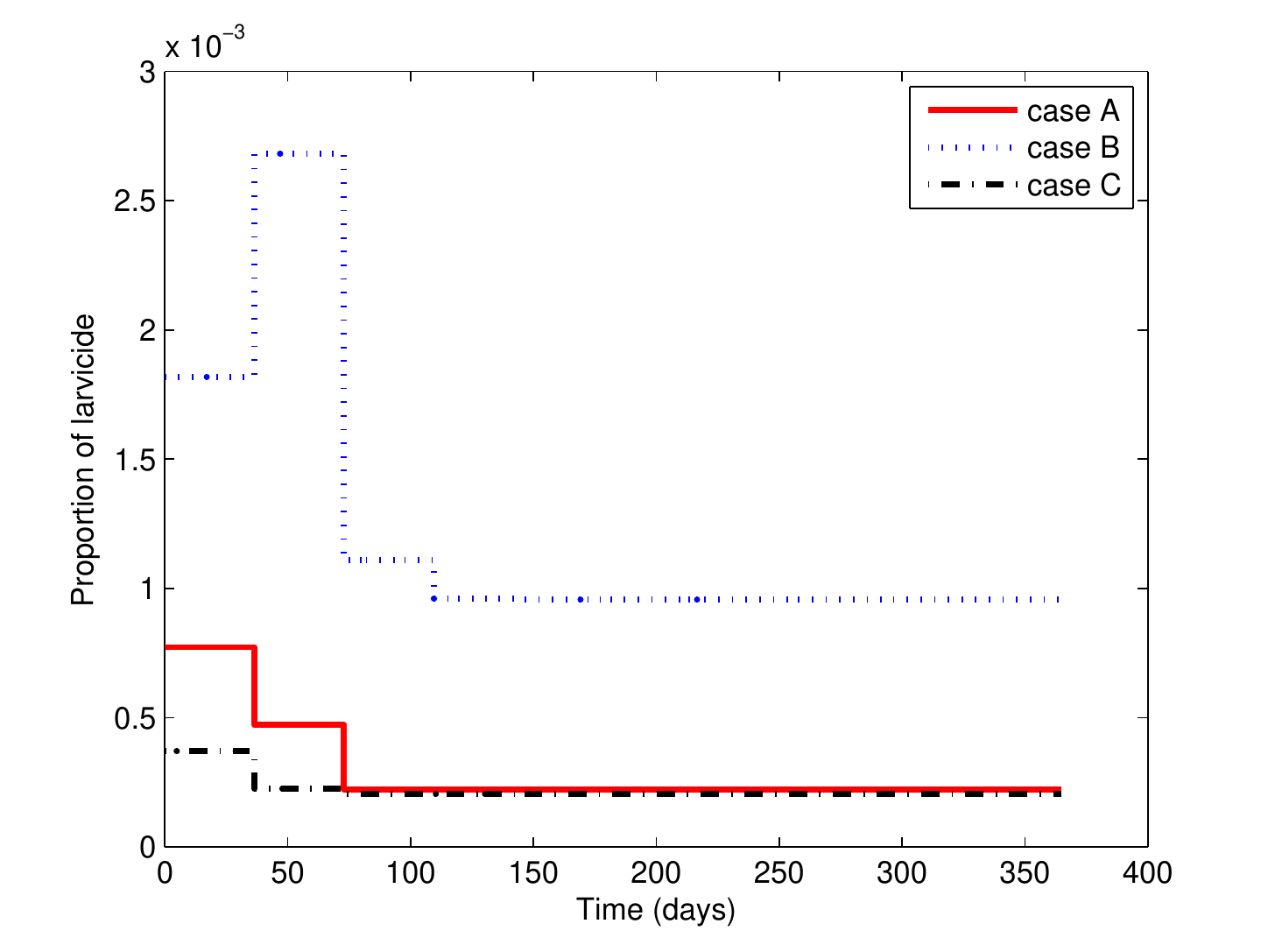}}\\
\subfigure[Mechanical control]{\label{cap6_proportion_mechanical_CasesABC}
\includegraphics[width=0.55\textwidth]{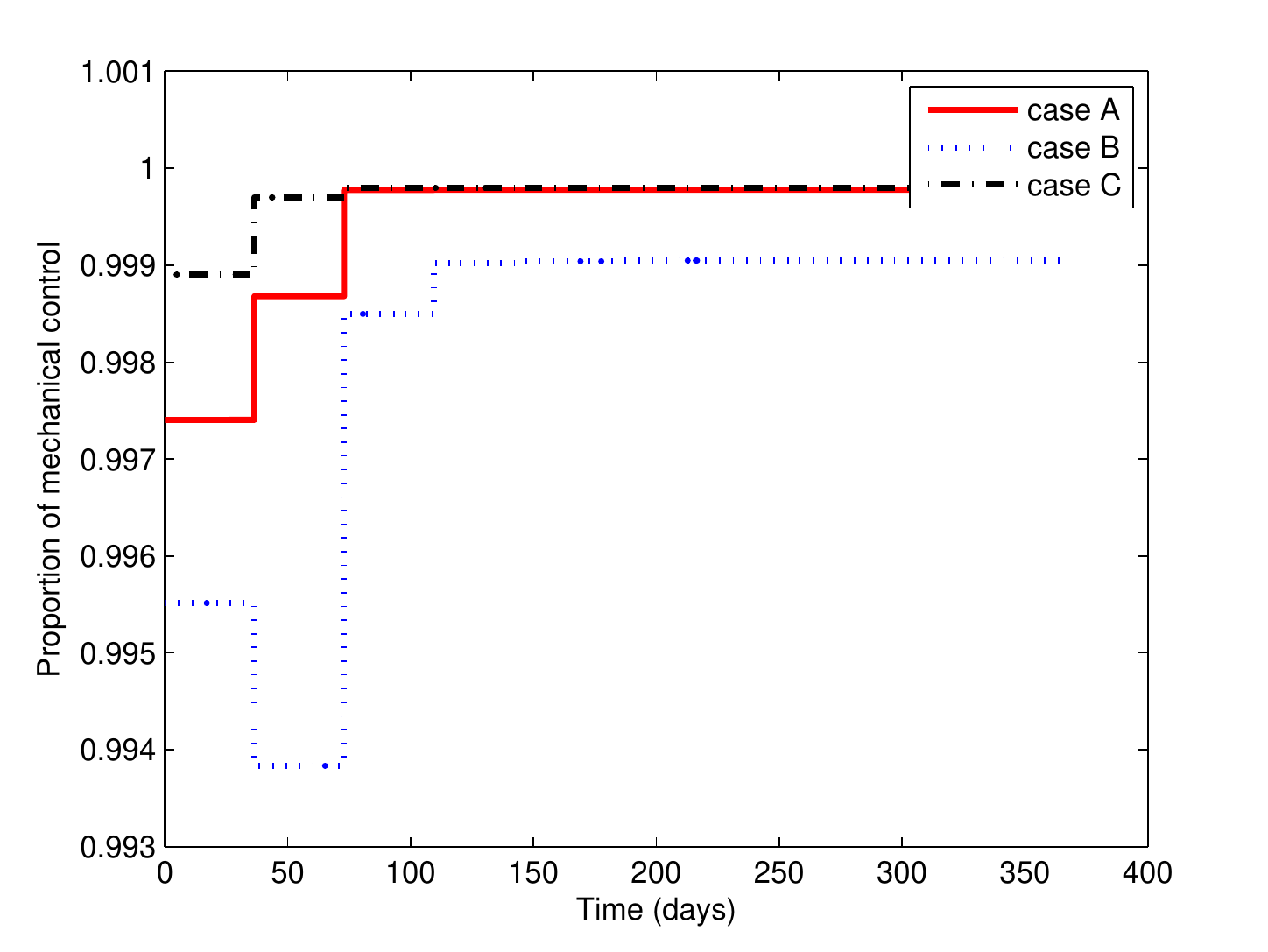}}
\caption{Proportion of control used in the three bioeconomic perspectives}
\label{cap6_all controls}
\end{figure}

Figure~\ref{cap6_human_infected_CasesABC} shows the number of infected human
in each bioeconomic perspective. We conclude that Case A and Case C are similar.
It can be explained by the low weight given to the cost of treatment (cases A and C)
when compared with the heavy weight given in case B. Figure~\ref{cap6_all controls}
presents the behavior of the controls for the A, B and C cases.
As adulticide is the control that has more influence in the model,
this is the one that most varies when the weights are changed.

\subsection{Only one control}

\begin{figure}
\centering
\subfigure[Optimal control]{\label{fcap6_adulticide_vs_all}
\includegraphics[width=0.48\textwidth]{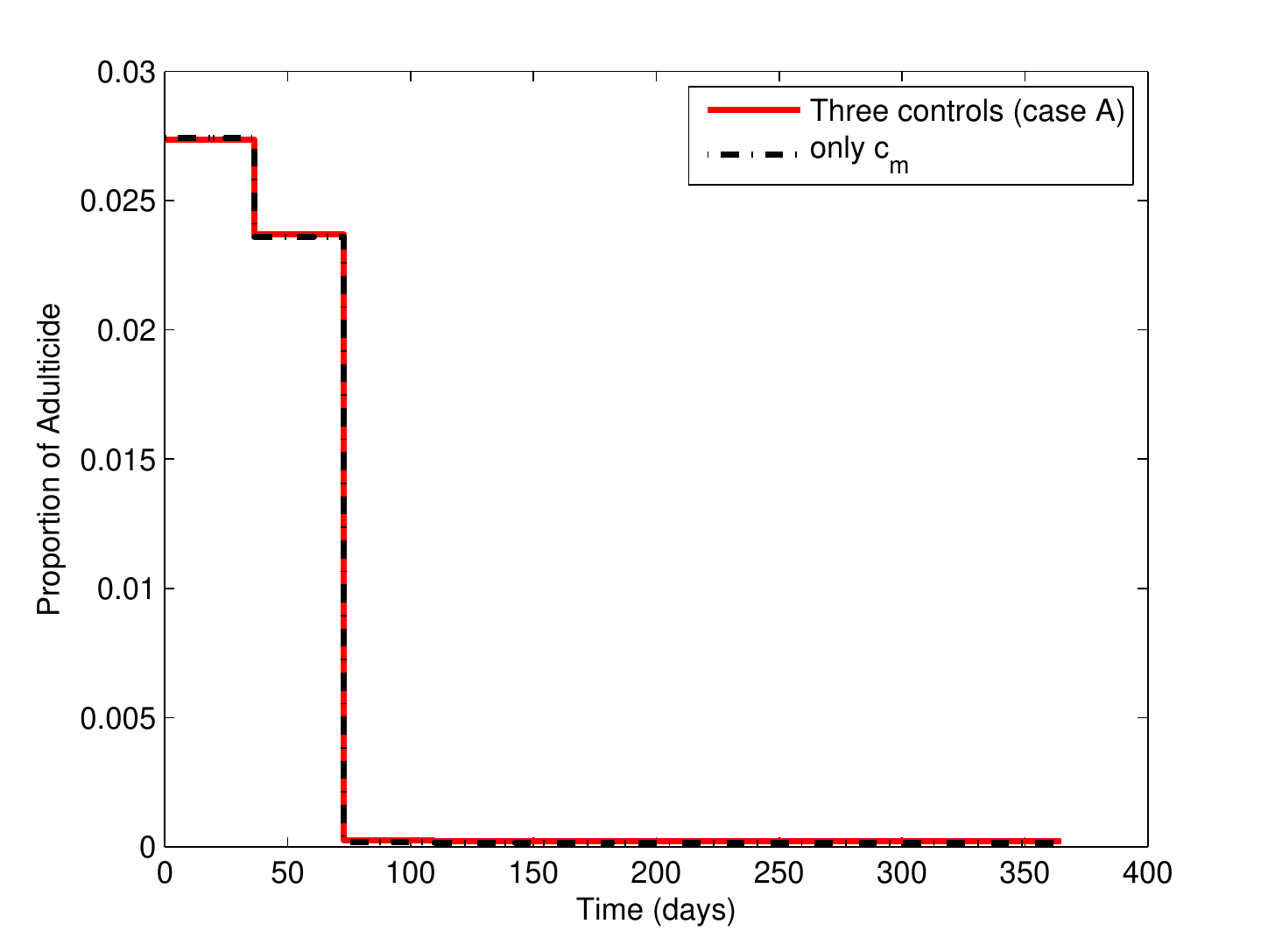}}
\subfigure[Proportion of infected human]{\label{cap6_infected_adulticide_vs_all}
\includegraphics[width=0.48\textwidth]{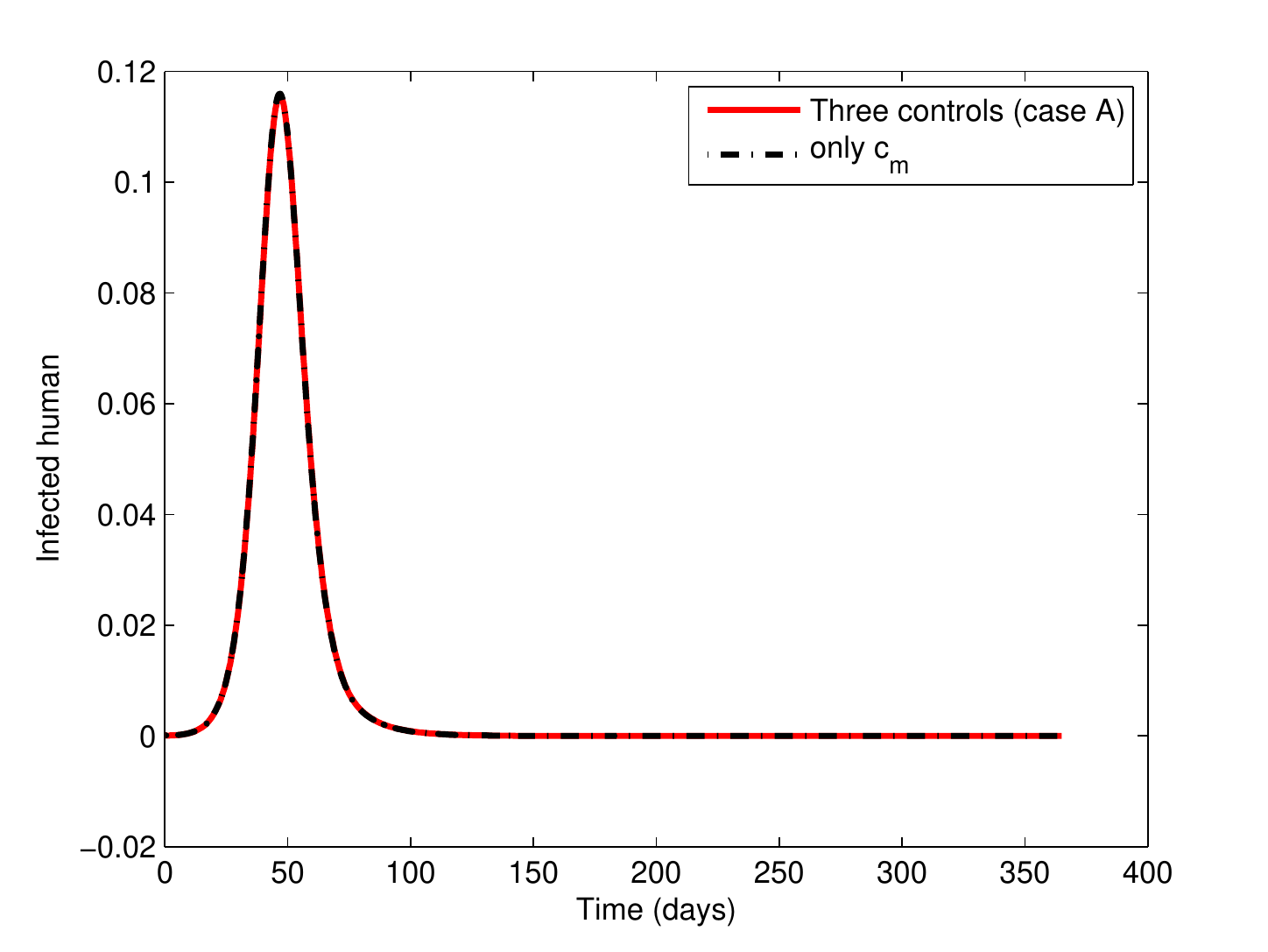}}
\caption{Optimal control and infected $I_h$ when all controls are considered
(solid line) and only adulticide control is taken into account (dashed line)}
\label{cap6_adulticide}
\end{figure}

\begin{figure}
\centering
\subfigure[Optimal control]{\label{cap6_larvicide_vs_all}
\includegraphics[width=0.48\textwidth]{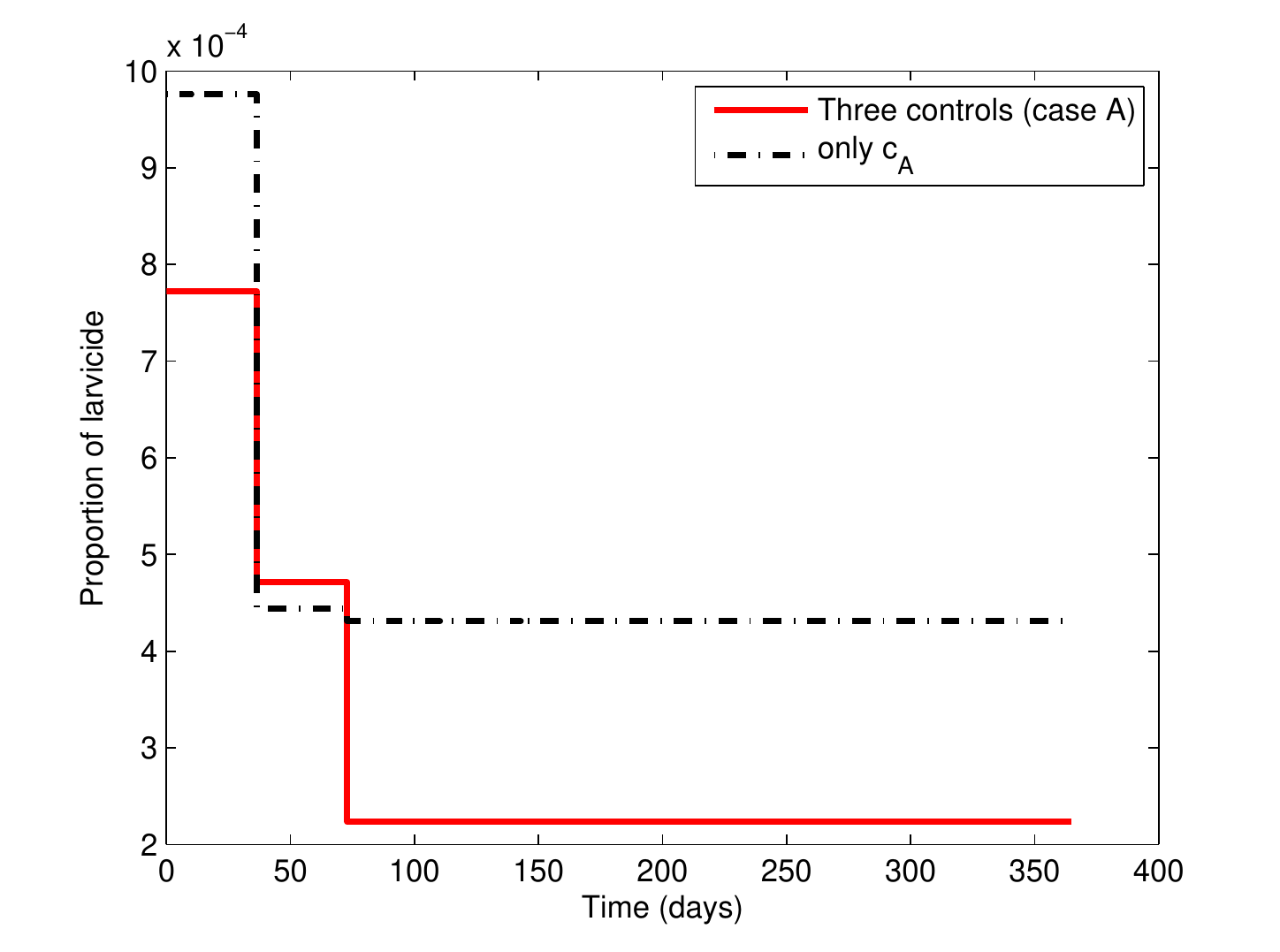}}
\subfigure[Proportion of infected human]{\label{cap6_infected_larvicide_vs_all}
\includegraphics[width=0.48\textwidth]{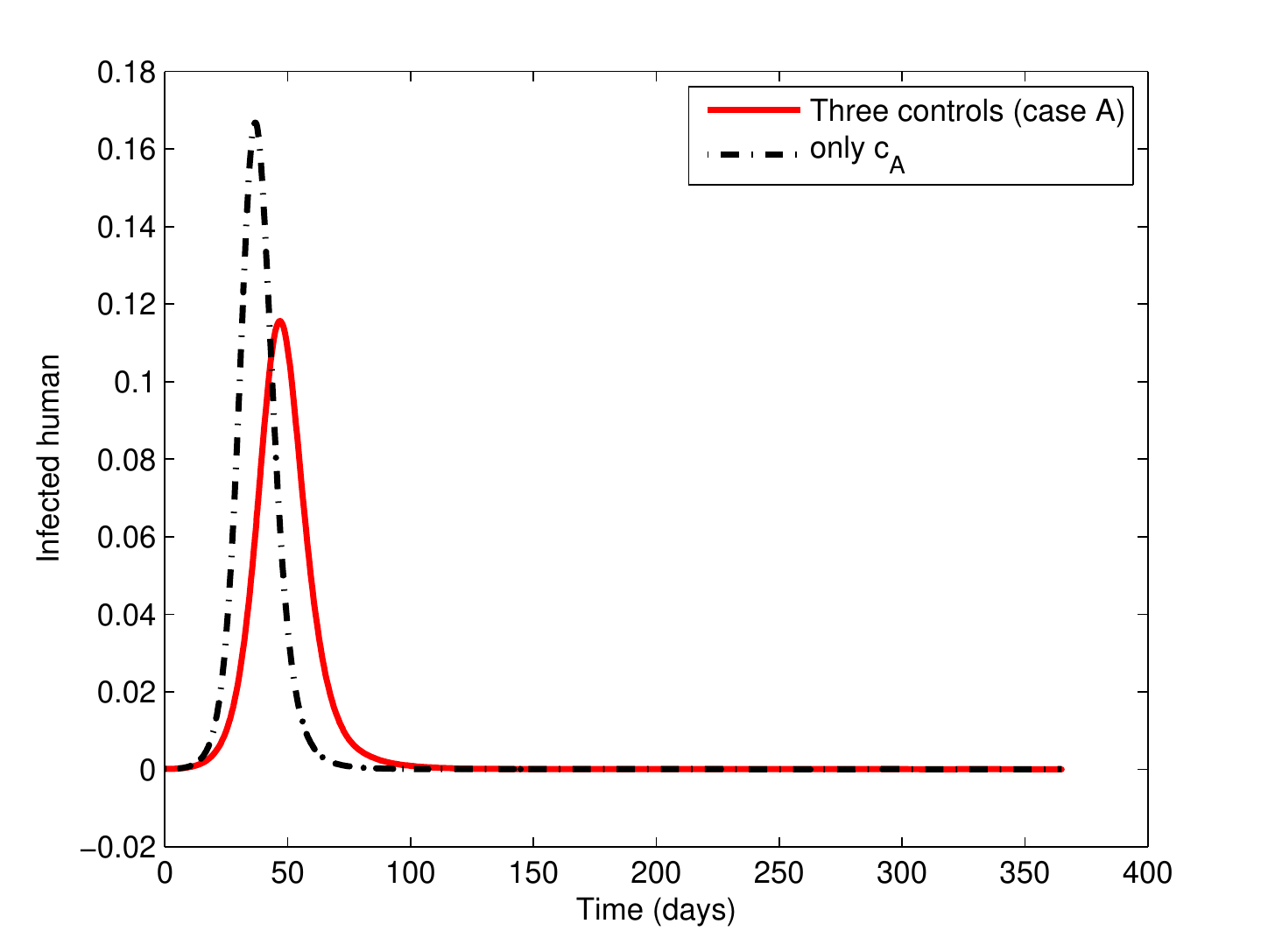}}
\caption{Optimal control and infected $I_h$ when all controls are considered
(solid line) and only larvicide control is taken into account (dashed line)}
\label{cap6_larvicide}
\end{figure}

\begin{figure}
\centering
\subfigure[Optimal control]{\label{cap6_mechanical_vs_all}
\includegraphics[width=0.48\textwidth]{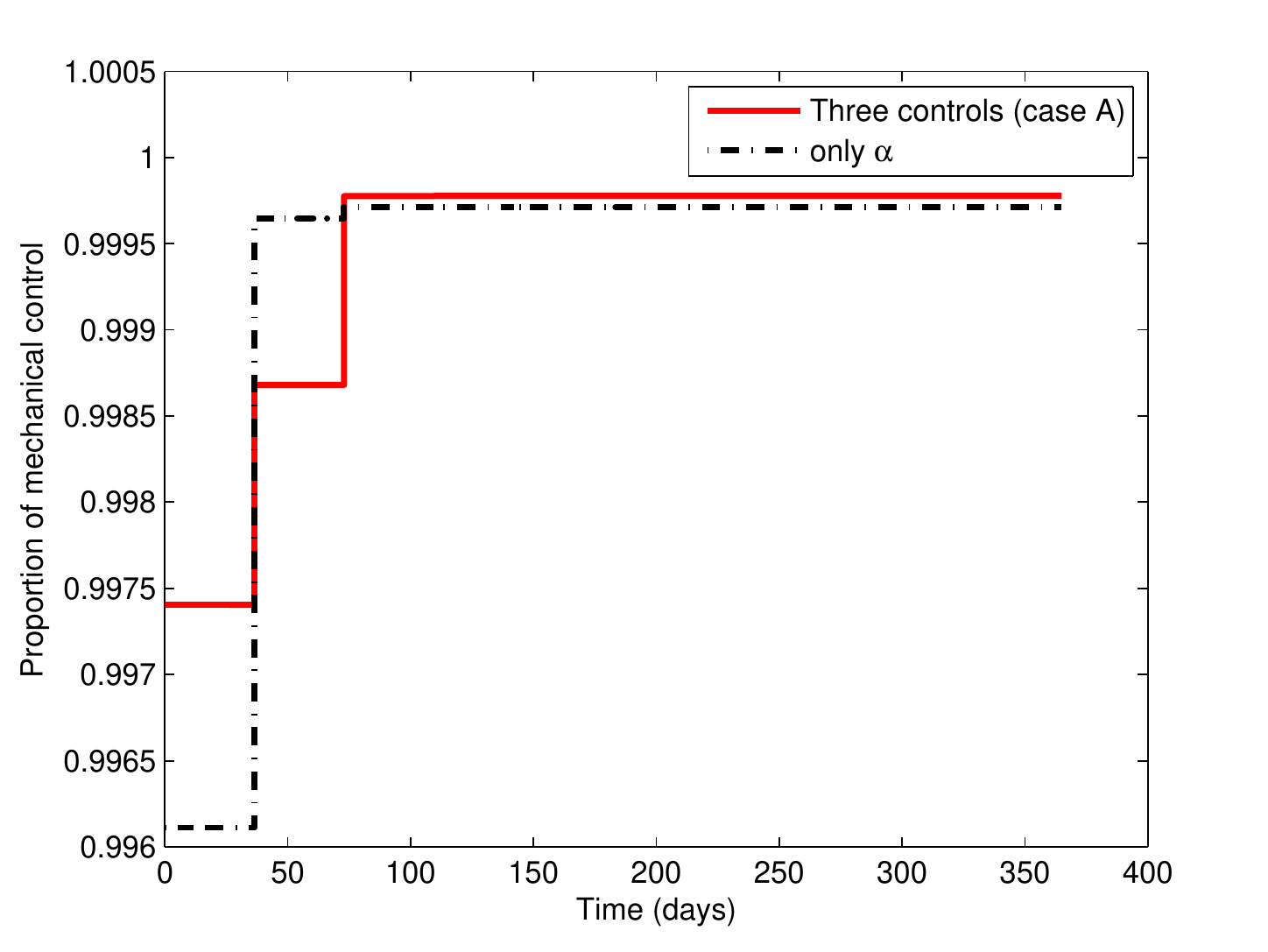}}
\subfigure[Proportion of infected human]{\label{cap6_infected_mechanical_vs_all}
\includegraphics[width=0.48\textwidth]{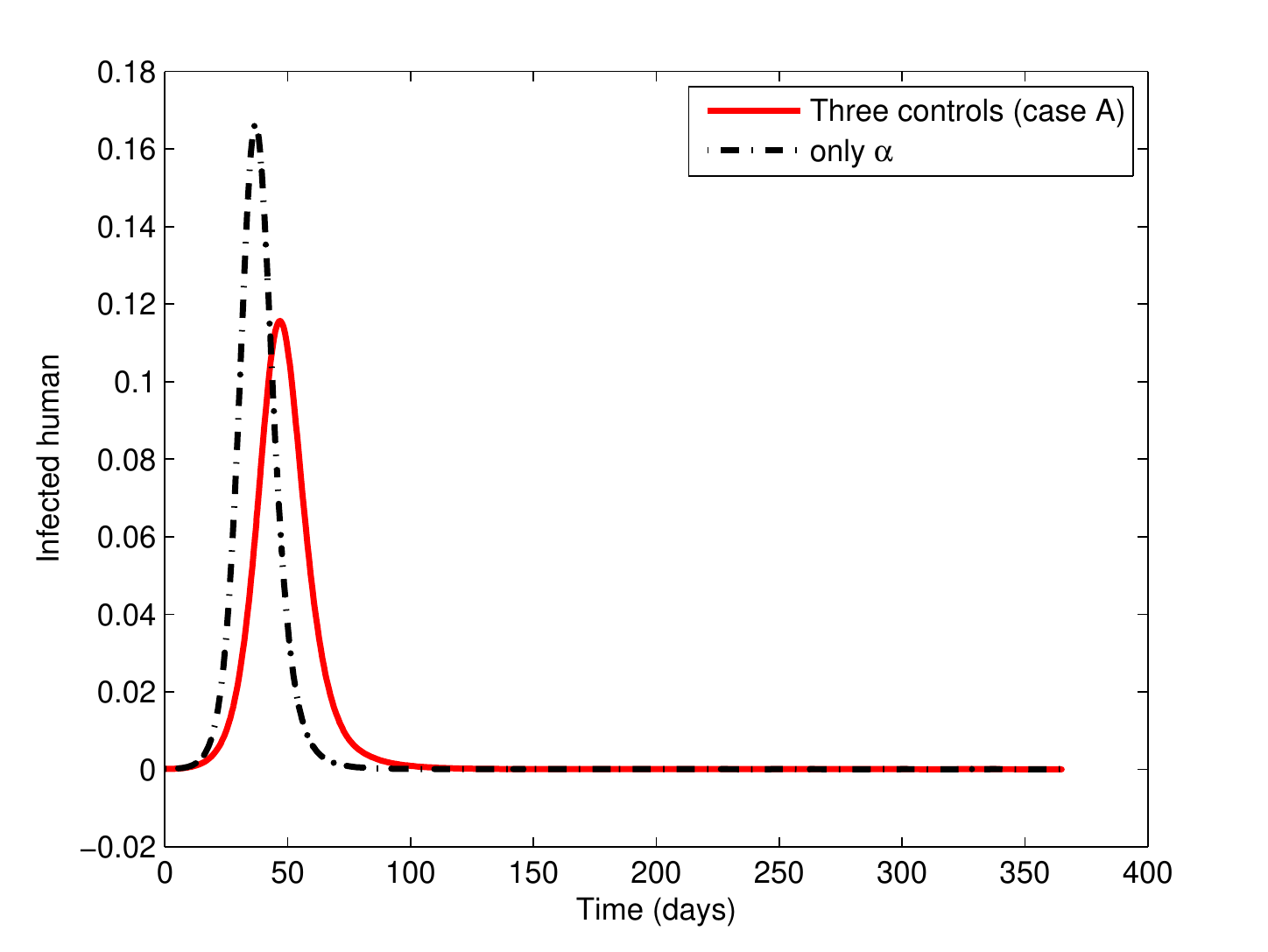}}
\caption{Optimal control and infected $I_h$ when all controls are considered
(solid line) and only mechanical control is taken into account (dashed line)}
\label{cap6_mechanical}
\end{figure}

A third strategy was tested: the functional was changed in order to study
the effects of each control when considered separately. Therefore,
the new functional also considers bioeconomic perspectives, but just
includes two variables: the costs with infected human (with $\gamma_D=0.5$)
and the costs with only one control (with $\gamma_i=0.5, i\in\{S,L,E\}$).
In Figure~\ref{cap6_adulticide} the proportion
of adulticide $(a)$ and infected humans $(b)$ are presented,
when the functional is $\int_{0}^{t_f}\left[\gamma_D i_h(t)^2
+\gamma_S c_m(t)^2\right]dt$.
Figures~\ref{cap6_larvicide} and \ref{cap6_mechanical} represent
the same simulations when the controls considered are larvicide
and mechanical control, respectively. It is possible to see that
the use of larvicide and mechanical control, used alone,
do not bring relevant influence to the control of the disease.


\section{Conclusions}
\label{sec:6}

Dengue disease breeds, even in the absence of fatal forms, significant
economic and social costs: absenteeism, debilitation and
medication. To observe and act at the onset of the epidemics,
can save lives and resources to governments.
Moreover, the under-reporting of dengue is probably
the most important barrier to obtaining an accurate assessment.

We presented a compartmental epidemiological
model for dengue, composed by a set of differential
equations. Simulations based on clean-up
campaigns to remove the vector breeding sites,
and also simulations on the application of insecticides
(larvicide and adulticide), were made. It was shown that even with
a low, although continuous, index of control over time,
the results are surprisingly positive. The adulticide
was the most effective control, from the fact that with
a low percentage of insecticide, the basic reproduction number
is kept below unit and the infected number of humans was smaller.

However, to bet only in adulticide is a risky decision. In some countries,
such as Mexico and Brazil, the prolonged use of adulticides has been increasing
the mosquito tolerance capacity to the product or even they become completely resistant.
In countries where dengue is a permanent threat, governments should
act with differentiated tools. It will be interesting to analyze these controls
in an endemic region and with several outbreaks.
We believe that the results will be quite different.
This is under investigation and will be addressed elsewhere.


\section*{Acknowledgements}

This work was partially supported by the Portuguese Foundation
for Science and Technology (FCT) through
the Ph.D. grant SFRH/BD/33384/2008 (Rodrigues)
and the R\&D units Algoritmi (Monteiro) and CIDMA (Torres).
Rodrigues and Torres were also supported
by the Strategic Project PEst-C/MAT/UI4106/2011.

The authors are very grateful to two anonymous referees 
for valuable remarks and comments, which significantly 
contributed to the quality of the paper.



\label{lastpage}

\end{document}